\documentclass[12pt,a4paper]{article}
\usepackage{bbm}
\usepackage{amsfonts}

\usepackage{amssymb}

\usepackage{amsmath}

\newtheorem{theorem}{Theorem}[section]

\newtheorem{lemma}[theorem]{Lemma}
\newtheorem{corollary}[theorem]{Corollary}

\begin{document}
\title{Gr\"{o}bner-Shirshov bases for braid groups in Adyan-Thurston generators\footnote{Supported by the
NNSF of China (No.10771077) and the NSF of Guangdong Province
(No.06025062).} }
\author{
Yuqun Chen and Chanyan Zhong  \\
\\
{\small \ School of Mathematical Sciences}\\
{\small \ South China Normal University}\\
{\small \ Guangzhou 510631}\\
{\small \ P. R. China}\\
{\small \ yqchen@scnu.edu.cn} \\
{\small \ chanyanzhong@yahoo.com.cn}}
\date{}
\maketitle \noindent\textbf{Abstract:} In this paper, we  give a
 Gr\"{o}bner-Shirshov basis of the braid group $B_{n+1}$ in Adyan-Thurston
 generators. We also deal with the braid group of type $\bf{B}_{n}$.
 As results, we obtain a new algorithm for getting the Adyan-Thurston normal form,
 and a new proof that the braid semigroup
$B^+_{n+1}$ is the subsemigroup in $B_{n+1}$.


\noindent {\bf AMS} Mathematics Subject Classification(2000):
 20F36, 20F05, 20F10,  16S15, 13P10

\section{Introduction}
Artin \cite{a26} invented a group $B_{n+1}$, the braid group on
$n+1$ strands
\begin{eqnarray*}
B_{n+1}=gp\langle \sigma_1,\dots,\sigma_{n}| \
\sigma_i\sigma_j=\sigma_j\sigma_i\ (i-1>j),
\sigma_{i+1}\sigma_i\sigma_{i+1}=\sigma_i\sigma_{i+1}\sigma_i
\rangle
\end{eqnarray*}
and solved the word problem for $B_{n+1}$. Markov \cite{ma45} and
Artin \cite{a47} found normal form for $B_{n+1}$ in Artin-Burau
generators
$$
s_{i,j}, \ s_{i,j}^{-1} \ (1\leq i<j\leq n), \ \sigma_{i}\ (1\leq
i\leq n),
$$
where $s_{i,j}=\sigma_{j-1}\cdots
\sigma_{i+1}\sigma_i^2\sigma_{i+1}^{-1}\cdots\sigma_{j-1}^{-1}$.
Markov and Artin gave another algorithm for the solution of the word
problem for $B_{n+1}$. Garside \cite{g69} found a normal form for
$B_{n+1}$ in Artin-Garside generators
$$
\Delta,\ \Delta^{-1},\ \sigma_i\ (1\leq i\leq n),
$$
where
$\Delta=\sigma_1\sigma_2\sigma_1\cdots\sigma_{n-1}\cdots\sigma_1\sigma_{n}\cdots\sigma_1$
and used the normal form for the positive solution of the conjugacy
problem for $B_{n+1}$. Birman-Ko-Lee \cite{bkl98} invented a new
presentation:
$$
B_{n+1}=gp\langle a_{ts}\ (1\leq s< t\leq n)|R\rangle,
$$
where $a_{ts}=(\sigma_{t-1}\cdots
\sigma_{s+1})\sigma_s(\sigma_{s+1}^{-1}\cdots\sigma_{t-1}^{-1})\
(1\leq s< t\leq n+1)$ and $R$ consists of the following relations
\[
\left \{
\begin{array}{ll}
a_{ts}a_{rq}=a_{rq}a_{ts}, & \mbox{for $(t-r)(t-q)(s-r)(s-q)>0$}, \\
a_{ts}a_{sr}=a_{tr}a_{ts}=a_{sr}a_{tr},  & \mbox{for $1\leq
r<s<t\leq n+1$}
 \end{array}
 \right.
\]
and found a normal form for $B_{n+1}$ in the new presentation. They
used the normal form for another algorithms for the solution of the
word and the conjugacy problems for $B_{n+1}$.

Bokut-Chainikov-Shum \cite{b007} found a Gr\"{o}bner-Shirshov
 basis for $B_{n+1}$ in Artin-Burau generators and as a corollary the
Markov-Artin normal form is followed. Bokut-Fong-Ke-Shiao \cite{bfk}
found a Gr\"{o}bner-Shirshov basis for the braid semigroup
\begin{eqnarray*}
B^+_{n+1}=sgp\langle \sigma_1,\dots,\sigma_{n}| \
\sigma_i\sigma_j=\sigma_j\sigma_i\ (i-1>j),
\sigma_{i+1}\sigma_i\sigma_{i+1}=\sigma_i\sigma_{i+1}\sigma_i
\rangle
\end{eqnarray*}
in the  Artin-Garside generators. Using this result, Bokut
\cite{bo08} found Gr\"{o}bner-Shirshov basis for $B_{n+1}$ in the
Artin-Garside generators. As a corollary, the Garside normal form
for $B_{n+1}$ is followed together with a new algorithm to reach the
Garside normal form of a braid. Bokut \cite{bo09} found a
Gr\"{o}bner-Shirshov basis for $B_{n+1}$ in the Birman-Ko-Lee
generators and hence a new algorithm and a new proof for
Birman-Ko-Lee normal form in $B_{n+1}$.

Braid group $B_n$ is a generalization of the symmetric group $S_n$,
which is the same as Artin group (a generalization of the Coxeter
group). The Coxeter graphs ${\bf{A}_n, \bf{B}_n, \bf{D}_n} $ (the
spherical type) are the same as Dynkin diagrams ${\bf{A}_n,
\bf{B}_n, \bf{D}_n}$ respectively. Hence, there are also finite
types ${\bf{A}_n, \bf{B}_n, \bf{D}_n, \bf{G}_2, \bf{F}_4,\bf{E}_6,
\bf{E}_7, \bf{E}_8}$ of braid groups.

The following preliminaries are related to the Gr\"{o}bner-Shirshov
bases for associative algebras.

Let $k$ be a field, $k\langle X\rangle$ the free associative algebra
over $k$ generated by $X$ and $ X^{*}$ the free monoid generated by
$X$, where the empty word is the identity which is denoted by 1. For
a word $w\in X^*$, we denote the length of $w$ by $|w|$. Let $X^*$
be a well ordered set. Let $f=\alpha \bar{f}+\Sigma\alpha_iu_i \in
k\langle X\rangle$, where $\alpha,\ \alpha_i\in k,\ \bar{f},\ u_i\in
X^*$ and $u_i<\bar{f}$. Then we call $\bar{f}$ the leading word and
$f$ monic if $\bar{f}$ has coefficient 1.

A well ordering $<$ on $X^*$ is  monomial if it is compatible with
the multiplication of words, that is, for $u, v\in X^*$, we have
$$
u < v \Rightarrow w_{1}uw_{2} < w_{1}vw_{2}\  \ for \  all \
 w_{1}, \ w_{2}\in  X^*.
$$
A standard example of monomial ordering on $X^*$ is the deg-lex
ordering to compare two words first by degree and then
lexicographically, where $X$ is a well ordered set.

Let $f$ and $g$ be two monic polynomials in \textmd{k}$\langle
X\rangle$ and $<$ a well ordering on $X^*$. Then there are two kinds
of compositions:

$ (i)$ If \ $w$ is a word such that $w=\bar{f}b=a\bar{g}$ for some
$a,b\in X^*$ with $|\bar{f}|+|\bar{g}|>|w|$, then the polynomial
 $ (f,g)_w=fb-ag$ is called the intersection composition of $f$ and
$g$ with respect to $w$.

$ (ii)$ If  $w=\bar{f}=a\bar{g}b$ for some $a,b\in X^*$, then the
polynomial $ (f,g)_w=f - agb$ is called the inclusion composition of
$f$ and $g$ with respect to $w$.

In $(f,g)_w$, $w$ is called the ambiguity of the composition.

Let $S\subset k\langle X\rangle$ such that every $s\in S$ is monic.
Then the composition $ (f,g)_w$ is called trivial modulo $ (S,w)$ if
$ (f,g)_w=\sum\alpha_i a_i s_i b_i$, where each $\alpha_i\in k$,
$a_i,b_i\in X^{*}, \ s_i\in S$ and $\overline{a_is_ib_i}<w$.

Generally, for $f,g\in k\langle X\rangle,\ f\equiv g \ \ \ mod(S,w)$
we mean $f-g=\sum\alpha_ia_is_ib_i,$ where every $\alpha_i\in k, \
s_i\in S,\ a_i,b_i\in X^*$ and $a_i\overline{s_i} b_i<w$.

$S$ is called a Gr\"{o}bner-Shirshov basis in $k\langle X\rangle$
with respect to the well ordering $<$ if any composition of
polynomials in $S$ is trivial modulo $S$.

 The following lemma was first proved
by Shirshov \cite{sh62} for  free Lie algebras (with deg-lex
ordering) (see also Bokut \cite{b72}). Bokut \cite{b76} specialized
the approach of Shirshov to associative algebras (see also Bergman
\cite{b}). For commutative polynomials, this lemma is known as
Buchberger's Theorem (see \cite{bu65, bu70}).

\ \

\noindent{\bf Composition-Diamond Lemma}  \ Let $k$ be a field, $A=k
\langle X|S\rangle=k\langle X\rangle/Id (S)$ and $<$ a monomial
ordering on $X^*$, where $Id (S)$ is the ideal of $k \langle
X\rangle$ generated by $S$. Then the following statements are
equivalent:
\begin{enumerate}
\item[ (i)] $S $ is a Gr\"{o}bner-Shirshov basis.
\item[ (ii)] $f\in Id (S)\Rightarrow \bar{f}=a\bar{s}b$
for some $s\in S$ and $a,b\in  X^*$.
\item[ (iii)] $Irr (S) = \{ u \in X^* |  u \neq a\bar{s}b ,s\in S,a ,b \in X^*\}$
is a $k$-basis of the algebra $A=k\langle X | S \rangle.$
\end{enumerate}

\ \

If a subset $S$ of $k\langle X \rangle$ is not a
Gr\"{o}bner-Shirshov basis then one can add to $S$ all nontrivial
compositions of polynomials of $S$ and continue this process
repeatedly in order to have a Gr\"{o}bner-Shirshov basis $S^{c}$
that contains $S$. Such a process is called the Shirshov algorithm.

Let $A=sgp\langle X|S\rangle$ be a semigroup presentation. Then $S$
is also a subset of $k\langle X \rangle$ and we can find
Gr\"{o}bner-Shirshov basis
 $S^{c}$, and $Irr(S^{c})$ is a normal form for $A$. We also call $S^{c}$ a
Gr\"{o}bner-Shirshov basis of $A$.

In this paper, we use the Composition-Diamond lemma to get the
Gr\"{o}bner-Shirshov normal form for the braid group $B_{n+1}$ in
Adyan-Thurston generators. It is exactly the left-greedy forms for
braid groups. We also use the same method to deal with the braid
group of type ${\bf \bf{B}_{n}}$.

\section{Gr\"{o}bner-Shirshov basis of the braid group $B_{n+1}$ in Adyan-Thurston
 generators}
In this section, we will give a Gr\"{o}bner-Shirshov basis of the
braid group $B_{n+1}$ in Adyan-Thurston
 generators.

Let $B_{n+1}$ denote the braid group of type  ${\bf \bf{A}_{n}}$.
Then
$$
B_{n+1}=gp\langle \sigma_1,\dots, \sigma_n\ |\
\sigma_j\sigma_i=\sigma_i\sigma_j\ (j-1>i),\
\sigma_{i+1}\sigma_i\sigma_{i+1}=\sigma_i\sigma_{i+1}\sigma_i\rangle.
$$

The symmetry group is as follow:
\begin{eqnarray*}
S_{n+1}=gp\langle s_1,\dots, s_n\ |\ s_i^2=1, s_js_i=s_is_j\
(j-1>i), s_{i+1}s_is_{i+1}=s_is_{i+1}s_i\rangle.
\end{eqnarray*}
Bokut and Shiao found the normal form for $S_{n+1}$ in the following
theorem.

\begin{theorem}(\cite{bs})
$N=\{s_{1i_1}s_{2i_2}\cdots s_{ni_n}|\ i_j\leq j+1\}$ is the
Gr\"{o}bner-Shirshov normal form for $S_{n+1}$ in generators
$s_i=(i, i+1)$ relative to the deg-lex ordering, where
$s_{ji}=s_js_{j-1}\cdots s_i\ (j\geq i),\ s_{jj+1}=1$.  \ \
$\square$
\end{theorem}

Let $\alpha\in S_{n+1}$ and $
\overline{\alpha}=s_{1i_1}s_{2i_2}\cdots s_{ni_n}\in N$ is the
normal form of $\alpha$. Define the length of $\alpha$ as
$|\overline{\alpha}|=l(s_{1i_1}s_{2i_2}\cdots s_{ni_n})$ and
$\alpha\perp\beta$ if
$|\overline{\alpha\beta}|=|\overline{\alpha}|+|\overline{\beta}|$.
Moreover, each $\overline{\alpha}\in N$ has a unique expression
$\overline{\alpha}=s_{_{l_1i_{l_1}}}s_{_{l_2i_{l_2}}}\cdots
s_{_{l_ti_{l_t}}}$, where each $s_{_{l_ji_{l_j}}}\neq1$. Such a $t$
is called the breath of $\alpha$.

We can easily get the following lemmas.

\begin{lemma}\label{l1}
Let $\alpha, \beta, \gamma\in S_{n+1}$. If
$|\overline{\alpha\beta\gamma}|=|\overline{\alpha}|+|\overline{\beta}|+|\overline{\gamma}|$,
then $\alpha \perp \beta \perp \gamma,\ \alpha \perp \beta\gamma$
and $\alpha\beta \perp \gamma$.  \ \ $\square$
\end{lemma}

\begin{lemma}\label{l2}
Let $\alpha, \beta, \gamma\in S_{n+1}$. If $\alpha\beta \perp
\gamma$ and $\alpha \perp \beta$, then $\alpha \perp \beta\gamma$
and $\beta \perp \gamma$.  \ \ $\square$
\end{lemma}

Now, we let
$$
B'_{n+1}=gp\langle r(\overline{\alpha}),\ \alpha\in
S_{n+1}\setminus\{1\}\ |\
r(\overline{\alpha})r(\overline{\beta})=r(\overline{\alpha\beta}),\
\alpha \perp \beta \rangle,
$$
where $r(\overline{\alpha})$ means a letter with the index
$\overline{\alpha}$.

Then $B_{n+1}\cong B'_{n+1}$. Indeed, define $\theta:
B_{n+1}\rightarrow B'_{n+1},\ \sigma_i\mapsto r(s_i)$ and $\theta':
B'_{n+1}\rightarrow B_{n+1},\ r(\overline{\alpha})\mapsto
\overline{\alpha}|_{s_i\mapsto \sigma_i}$. Then two mappings are
homomorphisms and $\theta\theta'=\mathbbm{l}_{B'_{n+1}},
\theta'\theta=\mathbbm{l}_{B_{n+1}}$. Hence,
$$
B_{n+1}=gp\langle r(\overline{\alpha}),\ \alpha\in
S_{n+1}\setminus\{1\}\ |\
r(\overline{\alpha})r(\overline{\beta})=r(\overline{\alpha\beta}),\
\alpha \perp \beta \rangle.
$$
Let $X=\{r(\overline{\alpha}),\ \alpha\in S_{n+1}\setminus\{1\}\}$.
The generator $X$ of $B_{n+1}$ is called Adyan-Thurston generator.
It is clear that each $r(\overline{\alpha})$ corresponds to a
positive braid which is non-repeating in Epstein at al's book
\cite{ep}.

Then the positive braid semigroup in generator $X$ is
$$
B_{n+1}^{+}=sgp\langle X\ |\
r(\overline{\alpha})r(\overline{\beta})=r(\overline{\alpha\beta}),\
\alpha \perp \beta \rangle.
$$

Let $s_1<s_2<\cdots<s_n$. Define
$r(\overline{\alpha})<r(\overline{\beta})$ if and only if
$|\overline{\alpha}|>|\overline{\beta}|$ or
$|\overline{\alpha}|=|\overline{\beta}|,\ \overline{\alpha}
<_{lex}\overline{\beta}$. It is clear that such an ordering on $X$
is well ordered. We will use the deg-lex ordering on $X^*$ in this
section.

\begin{theorem}\label{t2.4}
A Gr\"{o}bner-Shirshov basis of $B_{n+1}^{+}$ in Adyan-Thurston
generator $X$ relative to the deg-lex ordering on $X^*$ is:
\begin{eqnarray*}
r(\overline{\alpha})r(\overline{\beta})&=&r(\overline{\alpha\beta}),\ \ \ \alpha \perp \beta, \\
r(\overline{\alpha})r(\overline{\beta\gamma})&=&r(\overline{\alpha\beta})r(\overline{\gamma}),\
\ \ \alpha \perp \beta \perp \gamma.
\end{eqnarray*}
\end{theorem}

\textbf{Proof: } The composition of
$r(\overline{\alpha})r(\overline{\beta})$ and
$r(\overline{\beta})r(\overline{\gamma})$ would induce the relation
$r(\overline{\alpha})r(\overline{\beta\gamma})=r(\overline{\alpha\beta})r(\overline{\gamma})$
when $|\overline{\alpha\beta\gamma}|\neq
|\overline{\alpha}|+|\overline{\beta\gamma}|$.

 All possible ambiguities of compositions are:
\begin{enumerate}
\item[]\ $w_1=r(\overline{\alpha})r(\overline{\beta})r(\overline{\gamma}),\
\alpha\perp\beta\perp\gamma$,

\item[]\ $w_2=r(\overline{\alpha})r(\overline{\beta\gamma})r(\overline{\delta}),\
\alpha\perp\beta\perp\gamma,\ \beta\gamma\perp\delta$,

\item[]\ $w_3=r(\overline{\alpha})r(\overline{\beta})r(\overline{\gamma\delta}),\
\alpha\perp\beta\perp\gamma\perp\delta$,

\item[]\ $w_4=r(\overline{\alpha})r(\overline{\beta\gamma})r(\overline{\delta\mu}),\
\alpha\perp\beta\perp\gamma,\ \beta\gamma\perp\delta\perp\mu$.
\end{enumerate}
We only check the composition $w_4$. The others are similar.

Let
$f=r(\overline{\alpha})r(\overline{\beta\gamma})-r(\overline{\alpha\beta})r(\overline{\gamma}),\
g=r(\overline{\beta\gamma})r(\overline{\delta\mu})-r(\overline{\beta\gamma\delta})r(\overline{\mu})$.
Then, by Lemma \ref{l2}, $\gamma\perp\delta,\
\beta\perp\gamma\delta$ and
\begin{eqnarray*}
(f,g)_{w_4}&=&(r(\overline{\alpha})r(\overline{\beta\gamma})-r
(\overline{\alpha\beta})r(\overline{\gamma}))r(\overline{\delta\mu})
-r(\overline{\alpha})(r(\overline{\beta\gamma})r(\overline{\delta\mu})-
r(\overline{\beta\gamma\delta})r(\overline{\mu}))\\
&=&r(\overline{\alpha})r(\overline{\beta\gamma\delta})r(\overline{\mu})-
r(\overline{\alpha\beta})r(\overline{\gamma})r(\overline{\delta\mu})\\
&\equiv&r(\overline{\alpha\beta})r(\overline{\gamma\delta})r(\overline{\mu})-
r(\overline{\alpha\beta})r(\overline{\gamma\delta})r(\overline{\mu})\\
&\equiv&0.
\end{eqnarray*}

Hence the result holds.  \ \ $\square$

\ \

Let $\Delta=r(s_{11}s_{21}\cdots s_{n1})$. Then we have

\begin{lemma}(\cite{ep})
$r(s_i)\Delta=\Delta r(s_{n+1-i})$.  \ \ $\square$
\end{lemma}

In $B_{n+1}$, the following formulas hold.
\begin{enumerate}
\item[1)]\ $(\sigma_{i1}\sigma_{i+11}\cdots \sigma_{n1})(\sigma_{n+1-(i-1)n+1-(i-1)}
\sigma_{n+1-(i-2)n+1-(i-1)}\cdots \sigma_{nn+1-(i-1)})=\Delta$;

\item[2)]\ $(\sigma_{ii_1}\sigma_{jj_1})(\sigma_{j_1-11}\sigma_{i_12}\sigma_{i+22}
\sigma_{i+32}\cdots \sigma_{j2}s_{j+11}\cdots
\sigma_{n_1})=\sigma_{i1}\cdots \sigma_{n1}$;

\item[3)]\ $(\sigma_{ii_1}\sigma_{jj_1}\sigma_{kk_1})(\sigma_{k_1-11}\sigma_{j_12}\sigma_{i_1+13}
\sigma_{i+33}\sigma_{i+43}\cdots \sigma_{j+13}\sigma_{j+22}\cdots
\sigma_{k2}\sigma_{k+11}\cdots \sigma_{n_1})=\sigma_{i1}\cdots
\sigma_{n1}$.

\end{enumerate}

\begin{lemma}(\cite{ad84})
For any $\alpha\in S_{n+1}$, there exists an $E_{\alpha}\in S_{n+1}$
such that in $B_{n+1}$,
$r(\overline{\alpha})r(\overline{E_{\alpha}})=\Delta$.
\end{lemma}
\textbf{Proof: }  If $\alpha=s_i$, we set
$E_{\alpha}=s_{11}s_{21}\cdots s_{i-11}s_{i2}s_{i+11}\cdots s_{n1}$.

If $|\alpha|\geq 2$, we prove the result by induction on the breath
of $\alpha$.

By the above formulas and Lemma \ref{l1}, for any $\alpha\in
S_{n+1}\setminus\{1\}$, there exists $E_{\alpha}\in S_{n+1}$, such
that $\alpha E_{\alpha}=s_{11}s_{21}\cdots s_{n1}$ and
$|\overline{\alpha
E_{\alpha}}|=|\overline{\alpha}|+|\overline{E_{\alpha}}|=n(n+1)/2$.
Hence $r(\overline{\alpha})r(\overline{E_{\alpha}})=\Delta,\ \ \
\alpha \perp E_{\alpha}$.   \ \ $\square$

\ \

Now, we can represent the braid group as a semigroup:
$$
B_{n+1}=sgp\langle X,\  \Delta^{-1}\ |\
\Delta^{\varepsilon}\Delta^{-\varepsilon}=1,\ \varepsilon=\pm1,\
r(\overline{\alpha})r(\overline{\beta})=r(\overline{\alpha\beta}),\
\alpha \perp \beta \rangle.
$$

\begin{theorem}
A Gr\"{o}bner-Shirshov basis of $B_{n+1}$ in Adyan-Thurston
generator $X$ relative to the deg-lex ordering on $X^*$ is:
\begin{eqnarray*}
&&1)\ \ \ r(\overline{\alpha})r(\overline{\beta})=r(\overline{\alpha\beta}),\ \ \ \alpha \perp \beta, \\
&&2)\ \ \
r(\overline{\alpha})r(\overline{\beta\gamma})=r(\overline{\alpha\beta})r(\overline{\gamma}),\
\ \ \alpha \perp
\beta \perp \gamma,\\
&&3)\ \ \
r(\overline{\alpha})\Delta^{\varepsilon}=\Delta^{\varepsilon}r(\overline{\alpha}'),\
\ \
\overline{\alpha}'=\overline{\alpha}|_{s_i\mapsto s_{n+1-i}},\\
&&4)\ \ \ r(\overline{\alpha\beta})r(\overline{\gamma\mu})=\Delta
r(\overline{\alpha}')r(\overline{\mu}),\ \ \ \alpha \perp \beta
\perp \gamma \perp \mu,\ r(\overline{\beta\gamma})=\Delta,\\
&&5)\ \ \ \Delta^{\varepsilon}\Delta^{-\varepsilon}=1.
\end{eqnarray*}
\end{theorem}
\textbf{Proof: } We will prove that all possible compositions are
trivial modulo $S$. Denote by $(i\wedge j)_w$ the composition of the
type $i)$ and type $j)$ with respect to the ambiguity $w$. The
ambiguities $w$ of  all possible compositions are:
$$
\begin{array}{llll}
1\wedge 1\ \
r(\overline{\alpha})r(\overline{\beta})r(\overline{\gamma})&
1\wedge2\
r(\overline{\alpha})r(\overline{\beta})r(\overline{\gamma\mu})&
1\wedge 3\
r(\overline{\alpha})r(\overline{\beta})\Delta^\varepsilon& 1\wedge4\
r(\overline{\alpha})r(\overline{\beta\gamma})r(\overline{\mu\nu})\\
2\wedge1\
r(\overline{\alpha})r(\overline{\beta\gamma})r(\overline{\mu})&
2\wedge2\
r(\overline{\alpha})r(\overline{\beta\gamma})r(\overline{\mu\nu})&
2\wedge3 \
r(\overline{\alpha})r(\overline{\beta\gamma})\Delta^\varepsilon &
2\wedge4\
r(\overline{\alpha})r(\overline{\beta\gamma})r(\overline{\mu\nu})\\
3\wedge5\
r(\overline{\alpha})\Delta^{\varepsilon}\Delta^{-\varepsilon} &
4\wedge1\
r(\overline{\alpha\beta})r(\overline{\gamma\mu})r(\overline{\nu})&
4\wedge2\
r(\overline{\alpha\beta})r(\overline{\gamma\mu})r(\overline{\nu\omega})&
4\wedge3\
r(\overline{\alpha\beta})r(\overline{\gamma\mu})\Delta^\varepsilon \\
4\wedge4\
r(\overline{\alpha\beta})r(\overline{\gamma\mu})r(\overline{\nu\omega})&
5\wedge5\
\Delta^{\varepsilon}\Delta^{-\varepsilon}\Delta^{\varepsilon}
\end{array}
$$

We only check the composition $(4\wedge4)_w$. The others are
similar. Let $
f=r(\overline{\alpha\beta})r(\overline{\gamma\mu})-\Delta
r(\overline{\alpha}')r(\overline{\mu}),\
g=r(\overline{\gamma\mu})r(\overline{\nu\omega})-\Delta
r(\overline{\gamma}')r(\overline{\omega}), \
w=r(\overline{\alpha\beta})r(\overline{\gamma\mu})r(\overline{\nu\omega})$,
where $\alpha \perp \beta \perp \gamma \perp \mu \perp \nu \perp
\omega,\ r(\overline{\beta\gamma})=r(\overline{\mu\nu})=\Delta$.
Then
\begin{eqnarray*}
(f,g)_{w}&=&(r(\overline{\alpha\beta})r(\overline{\gamma\mu})-\Delta
r(\overline{\alpha}')r(\overline{\mu}))r(\overline{\nu\omega})
-r(\overline{\alpha\beta})(r(\overline{\gamma\mu})r(\overline{\nu\omega})-\Delta
r(\overline{\gamma}')r(\overline{\omega}))\\
&=&r(\overline{\alpha\beta})\Delta
r(\overline{\gamma}')r(\overline{\omega})-
\Delta r(\overline{\alpha}')r(\overline{\mu})r(\overline{\nu\omega})\\
&\equiv& \Delta
r(\overline{\alpha\beta}')r(\overline{\gamma}')r(\overline{\omega})-\Delta
r(\overline{\alpha }')\Delta r(\overline{\omega})\\
&\equiv& \Delta r(\overline{\alpha }')\Delta
r(\overline{\omega})-\Delta r(\overline{\alpha }')\Delta
r(\overline{\omega})\\
&\equiv&0.
\end{eqnarray*}

Hence the result holds.  \ \ $\square$

\begin{corollary}
Adyan-Thurston normal forms for $B_{n+1}$ are $\Delta^k
r(\overline{\alpha_1})\cdots r(\overline{\alpha_s})$, where $k\in
\mathbb{Z},\ r(\overline{\alpha_1})\cdots r(\overline{\alpha_s})$ is
minimal in deg-lex ordering.   \ \ $\square$
\end{corollary}

\noindent{\bf Remark}: Actually, the Adyan-Thurston normal forms for
the braid group are exactly the left greedy normal forms in Epstein
at al's book \cite{ep}.

\section{Gr\"{o}bner-Shirshov basis of the braid group of type ${\bf \bf{B}_{n}}$}

In this section, we will give a Gr\"{o}bner-Shirshov basis of the
braid group of type ${\bf{B}_{n}}$ by using the same method in
section 2.

Let $B(B_{n+1})$ denote the braid group of type ${\bf{B}_{n}}$. Then
\begin{eqnarray*}
B(B_{n+1})&=&gp\langle \sigma_1,\dots, \sigma_n\ |\
\sigma_j\sigma_i=\sigma_i\sigma_j\ (j-1>i),\
\sigma_{i+1}\sigma_i\sigma_{i+1}=\sigma_i\sigma_{i+1}\sigma_i\,\\
&&\ \ \ \ \ \
\sigma_n\sigma_{n-1}\sigma_n\sigma_{n-1}=\sigma_{n-1}\sigma_n\sigma_{n-1}\sigma_n\rangle.
\end{eqnarray*}
For the same as braid group of type ${\bf{A}_{n}}$, we define
\begin{eqnarray*}
G&=&gp\langle s_1,\dots, s_n\ |\ s_i^2=1, s_js_i=s_is_j\
(j-1>i),\ s_{i+1}s_is_{i+1}=s_is_{i+1}s_i\\
&& \ \ \ \ \ s_ns_{n-1}s_ns_{n-1}=s_{n-1}s_ns_{n-1}s_n\rangle.
\end{eqnarray*}
Then we can view $G$ as a semigroup with the same generators and
relations as group.

Let $s_1<s_2<\cdots<s_n$ and define the deg-lex ordering $<$ on
$S^*$, where $S=\{s_1,\dots,s_n\}$.

\begin{lemma}\label{l3.1}
A Gr\"{o}bner Shirshov basis of $G$ in generator $S$ relative to the
deg-lex ordering on $S^*$ is:
\begin{enumerate}
\item[1)]\ $s_i^2=1\ (1\leq i\leq n)$,
\item[2)]\ $s_js_i=s_is_j\ (j-1>i)$,
\item[3)]\ $s_{ji}s_j=s_{j-1}s_{ji}\ (1\leq i<j\leq n-1)$,
\item[4)]\ $s_{nj}s_{ni}=s_{n-1}s_{ni}s_{nj+1}\ (1\leq i\leq j\leq
n-1)$.
\end{enumerate}
\end{lemma}
\textbf{Proof: } We will prove that all possible compositions are
trivial modulo $S$. Denote by $(i\wedge j)_w$ the composition of the
type $i)$ and type $j)$ with respect to the ambiguity $w$. The
ambiguities $w$ of  all possible compositions are:
$$
\begin{array}{lllll}
1\wedge1\ \ s_{i}^3 & 1\wedge2\ s_j^2s_i & 1\wedge3\
s_js_{ji}s_j & 1\wedge4\ s_ns_{nj}s_{ni} & 2\wedge1\ s_js_i^2\\
2\wedge2\ s_ks_js_i & 2\wedge3 \ s_ks_{ji}s_j & 3\wedge1\
s_{ji}s_j^2 & 3\wedge2 \ s_{kj}s_ks_i& 3\wedge3\ s_{kj}s_{ki}s_k\\
4\wedge1\ s_{nj}s_{ni}s_i & 4\wedge2\ s_{nk}s_{nj}s_i & 4\wedge3\
s_{nk}s_{nj}s_i& 4\wedge4\ s_{nk}s_{nj}s_{ni}
\end{array}
$$

We only check the composition $(4\wedge4)_w$. The others are
similar. Let $w=s_{nk}s_{nj}s_{ni},\
f=s_{nk}s_{nj}-s_{n-1}s_{nj}s_{nk+1},\
g=s_{nj}s_{ni}-s_{n-1}s_{ni}s_{nj+1}$, where $1\leq i\leq j\leq
k\leq n-1$. Then
\begin{eqnarray*}
(f,g)_{w}&=&(s_{nk}s_{nj}-s_{n-1}s_{nj}s_{nk+1})s_{ni}
-s_{nk}(s_{nj}s_{ni}-s_{n-1}s_{ni}s_{nj+1})\\
&=&s_{nk}s_{n-1}s_{ni}s_{nj+1}-
s_{n-1}s_{nj}s_{nk+1}s_{ni}\\
&\equiv& s_{n-2}s_{nk}s_{ni}s_{nj+1}-
s_{n-1}s_{nj}s_{n-1}s_{ni}s_{nk+2}\\
&\equiv& s_{n-2}s_{n-1}s_{ni}s_{nk+1}s_{nj+1}-
s_{n-1}s_{n-2}s_{nj}s_{ni}s_{nk+2}\\
&\equiv& s_{n-2}s_{n-1}s_{ni}s_{n-1}s_{nj+1}s_{nk+2}-
s_{n-1}s_{n-2}s_{n-1}s_{ni}s_{nj+1}s_{nk+2}\\
&\equiv& s_{n-2}s_{n-1}s_{n-2}s_{ni}s_{nj+1}s_{nk+2}-
s_{n-2}s_{n-1}s_{n-2}s_{ni}s_{nj+1}s_{nk+2}\\
&\equiv&0.
\end{eqnarray*}

Hence the result holds.\ \ $\square$

\ \

 By using Lemma \ref{l3.1} and the Composition-Diamond
lemma, we have the following theorem.

\begin{theorem}
$N=\{s_{1i_1}s_{2i_2}\cdots s_{n-1i_{n-1}}s_{nj_1}\cdots s_{nj_k}|\
i_l\leq l+1,\ 1\leq j_1<j_2<\cdots<j_k\leq n,\ k\geq 0\}$ is the
Gr\"{o}bner-Shirshov normal form for $G$ in generator $S$ relative
to the deg-lex ordering on $S^*$, where $s_{ji}=s_js_{j-1}\cdots
s_i\ (j\geq i),\ s_{jj+1}=1$.  \ \ $\square$
\end{theorem}

\ \

Similar to the case of the braid group $B_{n+1}$ in the section 2,
we introduce the following notations.

Let $\alpha\in G$ and
$$\overline{\alpha}=s_{1i_1}s_{2i_2}\cdots s_{n-1i_{n-1}}s_{nj_1}\cdots
s_{nj_k}\in N
$$
is the normal form of $\alpha$. Define the length of $\alpha$ as
$|\overline{\alpha}|=l(s_{1i_1}s_{2i_2}\cdots
s_{n-1i_{n-1}}s_{nj_1}\cdots s_{nj_k})$ and $\alpha\perp\beta$ if
$|\overline{\alpha\beta}|=|\overline{\alpha}|+|\overline{\beta}|$.
Now, we let
$$
B(B'_{n+1})=gp\langle r(\overline{\alpha}),\ \alpha\in
G\setminus\{1\}\ |\
r(\overline{\alpha})r(\overline{\beta})=r(\overline{\alpha\beta}),\
\alpha \perp \beta \rangle.
$$

Then $B(B_{n+1})\cong B(B'_{n+1})$. Indeed, define $\theta:
B(B_{n+1})\rightarrow B(B'_{n+1}),\ \sigma_i\mapsto r(s_i)$ and
$\theta': B(B'_{n+1})\rightarrow B(B_{n+1}),\
r(\overline{\alpha})\mapsto \overline{\alpha}|_{s_i\mapsto
\sigma_i}$. Then two mappings are homomorphisms and
$\theta\theta'=\mathbbm{l}_{B(B'_{n+1})},
\theta'\theta=\mathbbm{l}_{B(B_{n+1})}$. Hence,
$$
B(B_{n+1})=gp\langle r(\overline{\alpha}),\ \alpha\in
G\setminus\{1\}\ |\
r(\overline{\alpha})r(\overline{\beta})=r(\overline{\alpha\beta}),\
\alpha \perp \beta \rangle.
$$

Let $X_1=\{r(\overline{\alpha}),\ \alpha\in G\setminus\{1\}\}$. Then
the positive braid semigroup of type ${\bf{B}_{n}}$ in generator
$X_1$ is:
$$
B(B_{n+1}^{+})=sgp\langle X_1\ |\
r(\overline{\alpha})r(\overline{\beta})=r(\overline{\alpha\beta}),\
\alpha \perp \beta \rangle.
$$

Define $r(\overline{\alpha})<r(\overline{\beta})$ if and only if
$|\overline{\alpha}|>|\overline{\beta}|$ or
$|\overline{\alpha}|=|\overline{\beta}|,\ \overline{\alpha}
<_{lex}\overline{\beta}$.

\ \

Similar to Theorem \ref{t2.4}, we have
\begin{theorem}
A Gr\"{o}bner-Shirshov basis of $B(B_{n+1}^{+})$ in generator $X_1$
relative to the deg-lex ordering on $X_1^*$ is:
\begin{eqnarray*}
r(\overline{\alpha})r(\overline{\beta})&=&r(\overline{\alpha\beta}),\ \ \ \alpha \perp \beta, \\
r(\overline{\alpha})r(\overline{\beta\gamma})&=&r(\overline{\alpha\beta})r(\overline{\gamma}),\
\ \ \alpha \perp \beta \perp \gamma.  \ \ \ \ \ \ \ \square
\end{eqnarray*}
\end{theorem}

Let $\Delta=r(s_{11}s_{21}\cdots s_{n-11}s_{n1}s_{n2}\cdots
s_{nn})$. Then we have

\begin{lemma}\label{l3}
$r(s_i)\Delta=\Delta r(s_{i})$.
\end{lemma}
\textbf{Proof: } We need only to show that in $B(B_{n+1})$
$$
\sigma_i(\sigma_{11}\sigma_{21}\cdots
\sigma_{n-11}\sigma_{n1}\sigma_{n2}\cdots
\sigma_{nn})=(\sigma_{11}\sigma_{21}\cdots
\sigma_{n-11}\sigma_{n1}\sigma_{n2}\cdots \sigma_{nn})\sigma_i.
$$

Suppose $i=n$. Then
\begin{eqnarray*}
&&\sigma_n(\sigma_{11}\sigma_{21}\cdots \sigma_{n-11}\sigma_{n1}\sigma_{n2}\cdots \sigma_{nn})\\
&=&(\sigma_{11}\sigma_{21}\cdots
\sigma_{n-21})\sigma_{n1}\sigma_{n1}\cdots
(\sigma_{n2}\cdots \sigma_{nn})\\
&=&(\sigma_{11}\sigma_{21}\cdots
\sigma_{n-21})\sigma_{n-1}\sigma_{n1}\sigma_{n2}(\sigma_{n2}
\sigma_{n3}\cdots \sigma_{nn})\\
&=&(\sigma_{11}\sigma_{21}\cdots
\sigma_{n-21})\sigma_{n-1}\sigma_{n-2}\sigma_{n1}
\sigma_{n2}\sigma_{n3}(\sigma_{n3}\cdots \sigma_{nn})\\
&=&(\sigma_{11}\sigma_{21}\cdots
\sigma_{n-21})\sigma_{n-1}\sigma_{n-2}\sigma_{n-3}\sigma_{n1}
\sigma_{n2}\sigma_{n3}\sigma_{n4}(\sigma_{n4}\cdots \sigma_{nn})\\
&=&\cdots\\
&=&(\sigma_{11}\sigma_{21}\cdots
\sigma_{n-21})\sigma_{n-1i}\sigma_{n1}\cdots \sigma_{nn-(i-1)}
(\sigma_{nn-(i-1)}\cdots \sigma_{nn})\\
&=&(\sigma_{11}\sigma_{21}\cdots
\sigma_{n-21})\sigma_{n-12}\sigma_{n1}\cdots \sigma_{nn-1}
(\sigma_{nn-1}\sigma_{nn})\\
&=&(\sigma_{11}\sigma_{21}\cdots
\sigma_{n-21})\sigma_{n-12}\sigma_{n1}\cdots
\sigma_{nn-2}\sigma_{n-1}
\sigma_{nn-1}(\sigma_{nn} \sigma_{nn})\\
&=&(\sigma_{11}\sigma_{21}\cdots
\sigma_{n-21}\sigma_{n-11}\sigma_{n1}\cdots \sigma_{nn-1}
\sigma_{nn})\sigma_{n}.
\end{eqnarray*}

Suppose $1\leq i\leq n-1$. Then
$$
\sigma_i(\sigma_{11}\sigma_{21}\cdots
\sigma_{n-11}\sigma_{n1}\sigma_{n2}\cdots \sigma_{nn})
=(\sigma_{11}\sigma_{21}\cdots
\sigma_{n-11})\sigma_{n-i}(\sigma_{n1}\sigma_{n2}\cdots
\sigma_{nn}).
$$
Since
\begin{eqnarray*}
&&(\sigma_{11}\sigma_{21}\cdots \sigma_{n-11}\sigma_{n1}\sigma_{n2}\cdots \sigma_{nn})\sigma_i\\
&=&(\sigma_{11}\sigma_{21}\cdots
\sigma_{n-11}\sigma_{n1}\sigma_{n2}\cdots
\sigma_{ni+1})\sigma_i(\sigma_{ni+2}\sigma_{ni+3}\cdots \sigma_{nn})\\
&=&(\sigma_{11}\sigma_{21}\cdots
\sigma_{n-11}\sigma_{n1}\sigma_{n2}\cdots
\sigma_{ni}\sigma_{ni})(\sigma_{ni+2}\sigma_{ni+3}\cdots \sigma_{nn})\\
&=&(\sigma_{11}\sigma_{21}\cdots
\sigma_{n-11}\sigma_{n1}\sigma_{n2}\cdots
\sigma_{ni-1})\sigma_{n-1}(\sigma_{ni}\sigma_{ni+1}\cdots \sigma_{nn})\\
&=&(\sigma_{11}\sigma_{21}\cdots
\sigma_{n-11}\sigma_{n1}\sigma_{n2}\cdots \sigma_{ni-2})
\sigma_{n-2}(\sigma_{ni-1}\sigma_{ni}\cdots \sigma_{nn})\\
&=&(\sigma_{11}\sigma_{21}\cdots
\sigma_{n-11}\sigma_{n1}\sigma_{n2}\cdots \sigma_{ni-3})
\sigma_{n-3}(\sigma_{ni-2}\sigma_{ni-1}\cdots \sigma_{nn})\\
&=&\cdots\\
&=&(\sigma_{11}\sigma_{21}\cdots
\sigma_{n-11}\sigma_{ni-(i-1)})\sigma_{n-(i-1)} (\sigma_{n2}
\sigma_{n3}\cdots \sigma_{nn})\\
&=&(\sigma_{11}\sigma_{21}\cdots \sigma_{n-11})\sigma_{n-i}
(\sigma_{n1}\sigma_{n2} \sigma_{n3}\cdots \sigma_{nn}),
\end{eqnarray*}
the result holds.   \ \ $\square$

\begin{lemma}
$r(s_i)r(E_i)=\Delta$, where
\begin{eqnarray*}
&&E_i=s_{11}s_{21}\cdots s_{i-11}s_{i2}s_{i+11}\cdots
s_{n-11}s_{n1}s_{n2}\cdots s_{nn},\ 1\leq i\leq n-1, \\
&&E_n=s_{11}s_{21}\cdots s_{n1}s_{n2}\cdots s_{nn-1}.
\end{eqnarray*}
\end{lemma}

\textbf{Proof: } By Lemma \ref{l3}, $s_n(s_{11}s_{21}\cdots
s_{n-11}s_{n1}s_{n2}\cdots s_{nn})=(s_{11}s_{21}\cdots
s_{n-11}s_{n1}s_{n2}\cdots s_{nn})s_n $ in $G$. Hence,
$s_nE_n=s_{11}s_{21}\cdots s_{n-11}s_{n1}s_{n2}\cdots s_{nn}$. For
$1\leq i\leq n-1$, since $ s_i(s_{11}s_{21}\cdots\\
s_{i-11}s_{i2}s_{i+11}\cdots s_{n-11})= s_{11}s_{21}\cdots s_{n-11}
$, $s_iE_i=s_{11}s_{21}\cdots s_{n-11}s_{n1}s_{n2}\cdots s_{nn}$.
But $|s_i|+|E_i|=|s_iE_i|$, we can get $r(s_i)r(E_i)=\Delta$. \ \
$\square$

\ \

Now, we can represent the braid group as a semigroup:
$$
B(B_{n+1})=sgp\langle X_1,\  \Delta^{-1}\ |\
\Delta^{\varepsilon}\Delta^{-\varepsilon}=1,\
r(\overline{\alpha})r(\overline{\beta})=r(\overline{\alpha\beta}),\
\alpha \perp \beta \rangle.
$$

Similar to the case of the braid group $B_{n+1}$ in the section 2,
we have the following theorem:
\begin{theorem}
A Gr\"{o}bner-Shirshov basis of $B(B_{n+1})$ in generator $X_1$
relative to the deg-lex ordering on $X_1^*$ is:
\begin{eqnarray*}
&&r(\overline{\alpha})r(\overline{\beta})=r(\overline{\alpha\beta}),\ \ \ \alpha \perp \beta, \\
&&r(\overline{\alpha})r(\overline{\beta\gamma})=r(\overline{\alpha\beta})r(\overline{\gamma}),\
\ \ \alpha \perp
\beta \perp \gamma,\\
&&r(\overline{\alpha})\Delta^{\varepsilon}=\Delta^{\varepsilon}r(\overline{\alpha}),\
\ \,\\
&&r(\overline{\alpha\beta})r(\overline{\gamma\mu})=\Delta
r(\overline{\alpha})r(\overline{\mu}),\ \ \ \alpha \perp \beta \perp
\gamma
\perp \mu,\ r(\overline{\beta\gamma})=\Delta,\  \overline{\alpha}=1 \mbox{ or }\ \overline{\mu}=1,\\
&&\Delta^{\varepsilon}\Delta^{-\varepsilon}=1.
\end{eqnarray*}
\end{theorem}

\begin{corollary}
The normal forms for $B(B_{n+1})$ are $\Delta^k
r(\overline{\alpha_1})\cdots r(\overline{\alpha_s}) (k\in
\mathbb{Z})$, where $r(\overline{\alpha_1})\cdots
r(\overline{\alpha_s})$ is minimal in deg-lex ordering.
\end{corollary}

\ \

\noindent{\bf Acknowledgement}: The authors would like to thank
Professor L.A. Bokut for his guidance, useful discussions and
enthusiastic encouragement in writing up this paper.

\end{document}